\documentclass[a4paper]{article}
\usepackage[latin1]{inputenc} % ou \usepackage[utf8]{inputenc}
\usepackage[T1]{fontenc} % ou \usepackage[OT1]{fontenc}
\usepackage{RR,RRthemes}
\usepackage{hyperref}
\usepackage{amssymb}
\usepackage[all]{xy}
\usepackage[latin1]{inputenc}
\usepackage{amsfonts}
\usepackage{amsmath}
\usepackage{amssymb}
\usepackage{url}
\usepackage{hyperref}%
\usepackage{graphicx}
\providecommand{\U}[1]{\protect\rule{.1in}{.1in}}
\usepackage[all]{xy}
\usepackage{color}
\RRdate{August 2010}
\newcommand{\point}{\mbox{\LARGE .}}
\def\blbx{\hbox{\vrule height 5pt width 5pt depth 0pt}\medskip}
\newcommand{\cqfd}{\hfill\blbx \\}
\newcommand{\proof}{\noindent\mbox{\bf Proof:}\\}

\def\calf{\mathcal{F}}
\def\E{\mathbb{E}}
\def\T{\mathcal{T}}

\def\PP{\mathbb{P}}
\def\LL{\mathbb{L}}
\def\QQ{\mathbb{Q}}
\def\EE{\mathbb{E}}

\newcommand{\Fa}{{\cal F}}

\def\Ea{\mathcal{E}}

\newtheorem{thm}{Theorem}[section]

\newtheorem{lem}[thm]{Lemma}

\newtheorem{prop}[thm]{Proposition}

\numberwithin{equation}{section}

\RRauthor{% les auteurs
Pierre Del Moral\thanks{ Centre INRIA Bordeaux et Sud-Ouest \& Institut de Math\'ematiques de Bordeaux , Universit\'e de Bordeaux I, 351 cours de la Lib\'eration 33405 Talence cedex, France, Pierre.Del-Moral@inria.fr}, 
 Peng Hu\thanks{ Centre INRIA Bordeaux et Sud-Ouest \& Institut de Math\'ematiques de Bordeaux , Universit\'e de Bordeaux I, 351 cours de la Lib\'eration 33405 Talence cedex, France, Peng.Hu@inria.fr}, Nadia Oudjane\thanks{EDF R \& D Clamart (nadia.oudjane@edf.fr)},
}
\authorhead{Del Moral, Hu and Oudjane}
\RRtitle{Enveloppe de Snell avec des crit\`eres d'optimalit\'e multiplicativement dependants de chemin}
\RRetitle{Snell envelope with path dependent multiplicative optimality criteria}
\titlehead{Snell envelope with path dependent multiplicative optimality criteria}
\RRresume{Nous analysons l'enveloppe de Snell avec des crit\`eres d'optimalit\'e multiplicativement dependants de chemin. Surtout pour ce cas, nous proposons une variation de backward r\'ecurrence de l'enveloppe de Snell qui permet d'\'etendre certains sch\'emas d'approximation classique \`a ce cas sp\'ecial. Dans ce cadre, nous proposons un sch\'ema d'approximation particule d'\'echantillonnage importance bas\'e sur un changement de mesure sp\'ecifique, destin\'e \`a concentrer l'effort de calcul dans les r\'egions soulign\'ees par les crit\`eres. Ce nouvel algorithme est th\'eoriquement \'etudi\'e. Nous fournissons des estimations non convergence asymptotique et de prouver que l'estimateur r\'esultant est surestim\'e.}
\RRabstract{
We analyze the Snell envelope with path dependent multiplicative optimality criteria. Especially for this case, we propose a variation of the Snell envelope backward recursion which allows to extend some classical approximation schemes to the multiplicatively path dependent case. In this framework, we propose an importance sampling particle approximation scheme based on a specific change of measure, designed to concentrate the computational effort in regions pointed out by the criteria. This new algorithm is theoritically studied. We provide non asymptotic convergence estimates and prove that the resulting estimator is high biased.   
}
\RRmotcle{enveloppe de Snell, option am\'ericain, mod\`ele particule, \'ev\'enements rares}
\RRkeyword{Snell envelope, american option, particle model, rare events}
\RRprojet{ALEA}
\RRdomaine{1} % cas du domaine numero 1
\RRthemeProj{cqfd} % theme du projet Apics
\RRdomaineProjBis{cqfd} % domaine du projet Apics
\RCBordeaux % centre de recherche Bordeaux - Sud Ouest
\begin{document}
 \RRNo{7360}

\makeRR   % cas d'un rapport de recherche

%% \makeRT % cas d'un rapport technique.
%% a partir d'ici, chacun fait comme il le souhaite
\section{Introduction}
%%%%%%%%%%%%%%%%%%%%%%

%The Multiplicative optimality criteria framework
%
The Snell envelope is related to the calculation of the optimal stopping time of  a random process based on a given optimality criteria. In this paper, we are interested in some complicated optimality criteria, especially the multiplicatively path dependent case. In other words, given a random process $(X_k)_{0\leq k\leq n}$ and some gain functions $(f_k)_{0\leq k\leq n}$ and $(G_k)_{0\leq k\leq n}$, we want to maximize the expected gain $\EE(f_\tau(X_\tau)\prod^{\tau-1}_{k=0}G_k(X_k))$ by choosing $\tau$ on a set of random stopping times $\T$.  For example, in finance, the multiplicative optimality criteria $(G_k)_{0\leq k\leq n}$ \footnote{
In present paper, if not specified, when one talks about the potential rare event $G_k$ or optimality criteria $G_k$, it means $G_k(X_k)$} can be interpreted as a discount factor related to a stochastic interest rate (taking then an exponential form), or as an obstacle for exotic options such as barriers in knock out options (taking then the form of indicator functions).
% $1_{A_k}(x_k)$ for a given zone $A_k$, at time $k$).

%Classical (non path dependent) framework
%
In the discrete time setting, these problems associated with Snell envelope are defined in terms of a given Markov process $(X_k)_{k\geq 0}$ taking values in some sequence of measurable state spaces  $(E_n,\mathcal{E}_k)_{k\geq 0}$ adapted to the natural filtration $\mathcal{F} = (\mathcal{F}_k)_{k\geq 0}$. We let $\eta_{0}=\mbox{\rm Law}(X_0)$ be the initial distribution on $E_{0}$,  and we denote by
$M_{k}(x_{k-1},dx_{k})$ the elementary Markov transition of the chain from $E_{k-1}$ into $E_{k}$. For a given time horizon $n$ and any $k\in \{0,\ldots, n\}$, we let $\T_k$ be the set of all stopping times $\tau$ taking values in $\{k,\ldots, n\}$. For a given sequence of non negative measurable functions $f_k$ on $E_k$, we define a target process $Z_k=f_k(X_k)$. 
Then $(U_k)_{0\leq k\leq n}$ the Snell envelope of process $(Z_k)_{0\leq k\leq n}$ is defined by a recursive formula:
$$
U_k = Z_k\vee\E(U_{k+1}|\calf_k)
$$
with terminal condition $U_n=Z_n$. The main property of the Snell envelope defined as above is 
$$
U_k = \sup_{\tau \in \T_k} \E(Z_\tau|\calf_k)=\E(Z_{\tau^*_k}|\calf_k)
~~\mbox{\rm with}
~~\tau^*_k = \min{\{k\leq j\leq n~:~ U_j=Z_j\}}\in \T_k
$$
Then the computation of the Snell envelope $(U_k)_{0\leq k\leq n}$ amounts to solve the following backward functional equation\footnote{
Consult the last paragraph of this section for a statement of the notation
used in this article.}
\begin{equation}\label{back}
u_k=f_k\vee M_{k+1}(u_{k+1}) 
\end{equation}
for any $0\leq k<n$ with the terminal condition $u_n=f_n$.

%How to analyse approximation schemes error ?
%
But at this level of generality, we can hardly have a closed solution of the function $u_k$. In this context, lots of numerical approximation schemes have been proposed. Most of them amounts to replace in recursion~(\ref{back}) the pair of functions and Markov transitions $(f_k,M_k)_{0\leq k\leq n}$ by some approximation model $(\widehat{f}_k,\widehat{M}_k)_{0\leq k\leq n}$ on some possibly reduced measurable subsets $\widehat{E}_k\subset E_k$.  
In paper~\cite{robust}, the authors provided a general robustness lemma to estimate the error related to the resulting approximation $\widehat{u}_k$ of the Snell eveloppe $u_k$, for several types of approximation models  $(\widehat{f}_k,\widehat{M}_k)_{0\leq k\leq n}$. 
\begin{lem}\label{lem:lemme1}
 For any $0\leq k<n$, on the state space $\widehat{E}_{k}$, we have that
$$
|u_{k}-\widehat{u}_{k}|  \le  \sum_{l= k}^{n} \widehat{M}_{k,l}
|f_l- \widehat{f}_l| + \sum_{l= k}^{n-1} \widehat{M}_{k,l}|(M_{l+1} -\widehat{M}_{l+1}) u_{l+1}|\ .
$$
\end{lem}
This lemma provides a natural way to compare and combine different approximation models. 
In the present paper, this Lemma will be applied in the specific framework of a multiplicative optimality criteria. 

%From the multiplicatively path dependant case to the standard case
%
Let us come back now to the multiplicatively path dependent case that we mentioned in the beginning of the article. Instead of $\EE(f_\tau(X_\tau))$ we want to maximize $\EE(f_\tau(X_\tau)\prod_{p=0}^{\tau-1} G_{p}(X_p))$ on the stopping times set $\T$. In this situation, a natural way is to consider the path $(X_0\dots X_k)_{0\leq k\leq n}$ as a new Markov chain $(\mathcal{X}_k)_{0\leq k\leq n}$ on path spaces and associate with transitions given for any
$\chi_{k-1}=(x_{0},\ldots,x_{k-1})\in (E_0\times \dots \times E_{k-1})$ and $\chi'_k=(x'_0,\ldots,x'_k)\in (E_0\times \dots \times E_{k})$ by the following formula
$$
\mathcal{M}_k(\chi_{k-1},d\chi'_k)
=\delta_{\chi_{k-1}}(d\chi'_{k-1})~M_k(x'_{k-1},dx'_k)\ .
$$
Then, let us denote by $\mathbf{u}_k(x_0\dots x_k)$, the Snell envelope defined with a path version of recursion~(\ref{back}):
\begin{equation}
\label{eq:RecursionPath}
\mathbf{u}_k(x_0,\dots,x_k)=\big[f_k(x_k)\prod^{k-1}_{p=0}G_p(x_p)\big]\vee \mathcal{M}_{k+1}(\mathbf{u}_{k+1})(x_0,\dots,x_k)\ ,
\end{equation}
for $0\leq k<n$ with terminal value $\mathbf{u}_n(x_0,\dots,x_n)=f_n(x_n)\prod^{n-1}_{p=0}G_p(x_p)$. 
%
%Numerical diificulties :Rare event modelization
%
At this stage, two difficulties may arise. First, the above recursion~(\ref{eq:RecursionPath}) seem to require the approximation of a $k+1$ dimensional function at each time step from $k=n-1$ up to $k=0$. Second, when the optimality criteria $G_p$ is localized in a specific region of $E_p$, for each $p$, then the product $\prod^{k-1}_{p=0}G_p(x_p)$ can be interpreted as a rare event.
% for which classical Monte Carlo methods are unaible to compute the expectation. 
Hence,  at first glance, the computation of Snell envelopes in the multiplicatively path dependent case seem to combine two additional numerical difficulties w.r.t. to the standard case, related to the computation of conditional expectations in a both high dimensional and rare event situation. 
%As one of the major problems in stochastic control and optimal stopping theory, the numerical solving of a normal Snell envelope requires already extensive calculations, most of Monte Carlo models explode in rare events situation.

These issues are considered in Section~2, of the present paper. The dimensionality problem is easily bypassed by considering an intermediate standard Snell envelope $(v_k)_{0\leq k\leq n}$, without path dependent criteria, which is directly related to the multiplicatively path dependent Snell envelope, by the relation $\mathbf{u}_k(x_0, x_1...x_k)=\prod_{p=0}^{k-1}G_p(x_p)v_k(x_k)$, for all $0\leq k\leq n$. Hence, computing the original Snell envelope $\mathbf{u}_k$  can be done by using one of the many approximation schemes developped for the standard (non path dependent) case. Then, to deal with the rare event problem, we propose a change of measure which allows to concentrate the computational effort in the regions of interest w.r.t. the criteria $(G_k)_{0\leq k\leq n-1}$. \\
In Section~3, we propose a Monte Carlo algorithm to compute the multiplicatively path dependent Snell envelope, on the base of this intermediate standard Snell envelope under a new equivalent measure defined in the previous section. 
This new approximation scheme is based on the stochastic mesh method introduced by M. Broadie and P. Glasserman in their seminal paper~\cite{glasserman} (see also~\cite{liuh}, for some recent refinements). The principal idea of original Broadie-Glasserman model is to make a change of probability, under the assumption that the Markov transitions $M_k(x,\cdot)$ are absolutely continuous w.r.t. some other measure $\eta_k$ on $E_k$, with positive Radon Nikodym derivatives $R_k(x,y)=\frac{dM_k(x,\point)}{d\eta_n}(y)$. But in most cases,  we do not know the density function of some good choice of $\eta_k$. So in~\cite{robust}, the authors  provide a variation of Broadie-Glasserman model that replaces not only $\eta_k$ but also the Radon-Nikodym derivatives $R_k$ with the approximation model $(\widehat{\eta}_k, \widehat{R}_k)$. The model introduced in the present article is an extension to multiplicatively path dependent functions. \\
In Section~4, the proposed Monte carlo algorithm is theorically analysed using an interacting particle system interpretation. We provide non asymptotic convergence estimates and prove that the resulting estimator is high biased.   

For the convenience of the reader, we end this introduction with some
notation used in the present article. 
We denote respectively by $\mathcal{P}
(E)$, and $\mathcal{B}(E)$, the set of all probability measures on some
measurable space $(E,\mathcal{E})$, and the Banach space of all bounded and
measurable functions $f$ equipped with the uniform norm $\Vert f\Vert $. 
We
let $\mu (f)=\int ~\mu (dx)~f(x)$, be the Lebesgue integral of a function $%
f\in \mathcal{B}(E)$, w.r.t. a measure $\mu \in \mathcal{P}(E)$. \\
We
recall that a bounded integral kernel $M(x,dy)$ from a measurable space $(E,%
\mathcal{E})$ into an auxiliary measurable space $(E^{\prime },\mathcal{E}%
^{\prime })$ is an operator $f\mapsto M(f)$ from $\mathcal{B}(E^{\prime })$
into $\mathcal{B}(E)$ such that the functions
\begin{equation*}
x\mapsto M(f)(x):=\int_{E^{\prime }}M(x,dy)f(y)
\end{equation*}%
are $\mathcal{E}$-measurable and bounded, for any $f\in \mathcal{B}%
(E^{\prime })$. In the above displayed formulae, $dy$ stands for an
infinitesimal neighborhood of a point $y$ in $E^{\prime }$. Sometimes, for indicator functions $f=1_A$, with $A\in\Ea$, we also use the 
notation $M(x,A):=M(1_A)(x)$.
The kernel $M$
also generates a dual operator $\mu \mapsto \mu M$ from $\mathcal{M}(E)$
into $\mathcal{M}(E^{\prime })$ defined by $(\mu M)(f):=\mu (M(f))$. A
Markov kernel is a positive and bounded integral operator $M$ with $M(1)=1$.
Given a pair of bounded integral operators $(M_{1},M_{2})$, we let $%
(M_{1}M_{2})$ be the composition operator defined by $%
(M_{1}M_{2})(f)=M_{1}(M_{2}(f))$. 
Given a sequence of  bounded integral operators $M_{n}$  from some state space 
$E_{n-1}$ into another $E_{n}$, we set $M_{k,l}:=M_{k+1}M_{k+2}\cdots M_l $,
for any $k\leq l$, with the convention $M_{k,k}=Id$, the identity  operator.
 In the context of finite state spaces, these integral operations coincide with the traditional matrix operations on multidimensional state spaces.\\
We also assume
that the reference Markov  chain $X_n$ with  initial distribution  $\eta_{0}\in\mathcal{P}(E_0)$, 
and elementary transitions 
$M_{n}(x_{n-1},dx_{n})$ from $E_{n-1}$ into $E_n$ is defined on some filtered probability space
$(\Omega,\calf,\PP_{\eta_0})$, and we
use the notation $\EE_{\PP_{\eta_0}}$ to denote the expectations
w.r.t. $\PP_{\eta_0}$. In this notation, 
for all $n\geq 1$ and for any $f_n\in \mathcal{B}(E_n)$, we have that
$$
\E_{\PP_{\eta_0}}\left\{ f_n(X_{n})|\calf_{n-1}\right\} = M_n f_n(X_{n-1}):=\int_{E_n}~M_n(X_{n-1},dx_n)~f_n(x_n)
$$
with the $\sigma$-field  $\Fa_n=\sigma(X_0,\ldots,X_n)$ generated by the sequence of random
variables $X_p$, from the origin $p=0$ up to the time $p=n$.  We also use the conventions $\prod_{\emptyset}=1$, and $\sum_{\emptyset}=0$.

\section{Snell envelope with multiplicatively path dependent functions and change of measure}\label{sub1}
%
%\section{Analysis of Snell envelope on multiplicatively path dependent functions}\label{sub1}
Suppose $(X_k)_{0\leq k\leq n}$ is a Markov chain on continuous state spaces $(E_k, \mathcal{E}_{k})_{0\leq k\leq n}$ with an initial distribution $\eta_{0}$ on $E_{0}$, a collection of Markov transitions $M_{k}(x_{k-1},$ $dx_{k})$ from $E_{k-1}$ to $E_k$ and a given final time horizon $n$. We also assume that the chain $X_k$ is defined on a filtered probability space $(\Omega ,\mathcal{F},\PP)$. 
 In this situation, the historical process $\mathcal{X}_k:=(X_0,\dots,X_k)$ can be seen as a Markov chain with transitions given for any
$\chi_{k-1}=(x_{0},\ldots,x_{k-1})\in E_0\times \dots \times E_{k-1}$ and $\chi'_k=(x'_0,\ldots,x'_k)\in E_0\times \dots \times E_{k}$ by the following formula
$$
\mathcal{M}_k(\chi_{k-1},d\chi'_k)
=\delta_{\chi_{k-1}}(d\chi'_{k-1})~M_k(x'_{k-1},dx'_k)\ .
$$
We denote by $(\PP_k)_{0\leq k\leq n}$ a sequence of probabilities of path $(\mathcal{X}_k)_{0\leq k\leq n}$.
For a given collection of real valued functions $(f_k)_{0\leq k\leq n}$  and $(G_k)_{0\leq k\leq n}$, defined on $(E_k)_{0\leq k\leq n}$, we define a class of real valued functions $(F_k)_{0\leq k\leq n}$ defined on the product spaces $(E_0\times \dots \times E_{k})_{0\leq k\leq n}$ by 
 $$
 F_k(x_0, \cdots ,x_k):=f_k(x_k)\prod_{0\leq p\leq k-1}G_p(x_p)\ ,\quad \textrm{for all}\   0\leq k\leq n\ .
 $$ 
To maximize the expected gain $\EE(F_\tau(\mathcal{X}_\tau))$ w.r.t. $\tau$ in a set of random stopping times $\T$, one is interested in computing the Snell envelope $(\mathbf{u}_k)_{0\leq k\leq n}$ associated to the gain functions $(F_k)_{0<k\leq n}$ and solution to the following recursion
\begin{equation}
\label{eq:u}
\left\{
\begin{array}{ll}
\mathbf{u}_n(x_0, \cdots ,x_n)=F_n(x_0, \cdots, x_n)\\
\mathbf{u}_{k}(x_0, \cdots , x_{k})=F_{k}(x_0, \cdots , x_{k})\vee \mathcal{M}_{k+1}(\mathbf{u}_{k+1})(x_0,\dots ,x_{k}), \forall ~0\leq k\leq n-1
\end{array}
\right .
\end{equation}
Now, let us consider the standard (non path dependent) Snell envelope $(v_k)_{0\leq k\leq n}$ associated to the gain functions $(f_k)_{0\leq k\leq n}$ and satisfying the following recursion
\begin{equation}\label{RSE0}
\left\{
\begin{array}{lll}
v_n(x_n)&=&f_n(x_n)\\
v_k(x_k)&=&f_k(x_k)\vee\big[ G_k(x_k)  M_{k+1}(v_{k+1})(x_{k})\big]\ ,\ \textrm{for all}\ 0\leq k\leq n-1\ .
\end{array}
\right.
\end{equation}
For all $0\leq k\leq n$, let us denote by $\mathbf{v}_k$ the real valued functions defined on $E_0\times \dots \times E_{k}$, such that $\mathbf{v}_k(x_0, \cdots , x_k):=v_k(x_k)\prod_{p=0}^{k-1}G_p(x_p)$. 
By construction, one can easily check  that for all $0\leq k\leq n$, 
$\mathbf{u}_k\equiv \mathbf{v}_k $
and in particular $\mathbf{u}_0(x_0)=v_0(x_0)$. 
Indeed, one can verify that $(\mathbf{v}_k)_{0\leq k\leq n}$ follows the same recursion~(\ref{eq:u}) as $(\mathbf{u}_k)_{0\leq k\leq n}$. First, we note that they share the same terminal condition, 
\begin{eqnarray*}
\mathbf{v}_n(x_0, \cdots , x_n)&=&v_n(x_n)\prod_{p=0}^{n-1}G_p(x_p)=f_n(x_k)\prod_{p=0}^{n-1}G_p(x_p)\\
&=&F_n(x_0, \cdots, x_n)=\mathbf{u}_n(x_0, \cdots , x_n)\ .
\end{eqnarray*}
Then at time step $k$, we observe that they follow the same recursion
\begin{eqnarray*}
&&\mathbf{v}_k(x_0, \cdots , x_k)\\
&=&v_k(x_k)\prod_{p=0}^{k-1}G_p(x_p)\\
&=&\Big [f_k(x_k)\prod_{p=0}^{k-1}G_p(x_p)\Big]\vee\Big[ \int M_{k+1}(x_k,dx_{k+1})v_{k+1}(x_{k+1})G_k(x_k) \prod_{p=0}^{k-1}G_p(x_p)\Big]\\
%&=&f_k(x_k)\big(\prod_{p=0}^{k-1}G_p(x_p)\big)\vee\left( G_k(x_k) \big(\prod_{p=0}^{k-1}G_p(x_p)\big)\int M_{k+1}(x_k,dx_{k+1})v_{k+1}(x_{k+1})\right)
&=&F_{k}(x_0, \cdots , x_{k})\vee \mathcal{M}_{k+1}(\mathbf{v}_{k+1})(x_0,\dots ,x_{k})\ .
\end{eqnarray*}
Now that we have underlined the link between $\mathbf{u}_k$ and $v_k$, our aim is then to compute the latter.
%Now if the criteria $G_k$ are localized in a specific region of $E_k$,  the recursion~(\ref{RSE0}) implies that 
%one concludes that the conditional expectation $M_{k+1}(v_{k+1})(x_k)$
The recursion~(\ref{RSE0}) implies that it is not relevant to compute precisely the conditional expectation $M_{k+1}(v_{k+1})(x_k)$ when the value of the criteria $G_k(x_k)$ is zero or very small.  Similarly, notice that $v_{k+1}$ is likely to reach high values when $G_{k+1}$ does, hence  from a variance reduction point of view, when approximating the conditional expectation $M_{k+1}(v_{k+1})(x_k)$ by a Monte Carlo method, it seems relevant to concentrate the simulations in the regions of $E_{k+1}$ where $G_{k+1}$ reaches high values. 
%
%
%As one of the major problems in stochastic control and optimal stopping theory, the numerical solving of a normal Snell envelope requires already extensive calculations, most of Monte Carlo models explode in rare events situation.
Hence, to avoid the potential rare events $G$, we propose to consider the following change of measure one the measurable product space $(E_0\times\cdots \times E_n, \mathcal{E}_0\times\cdots \times \mathcal{E}_n)$, 
\begin{equation}\label{chgpro}
d\QQ _n=\frac{1}{Z_n}\left[\prod_{k=0}^{n-1}G_k\right]d\PP _n\ ,\quad \textrm{with}\quad 
%\end{equation}
%\begin{equation}\label{chg}
Z_n=\E\left(\prod_{k=0}^{n-1}G_k(X_k)\right)=\prod_{k=0}^{n-1}\eta_k(G_k)\ ,
\end{equation}
%\begin{equation}\label{chgpro}
%\QQ _n(d(x_0\ldots x_n))=\frac{1}{Z_n}\left[\prod_{k=0}^{n-1}G_k(x_k)\right]\PP _n(d(x_0\ldots x_n))
%\end{equation}
% \begin{equation}\label{chg}
% Z_n=\E\left(\prod_{k=0}^{n-1}G_k\right)=\prod_{k=0}^{n-1}\eta_k(G_k)\ ,
% \end{equation}
where $\eta_k$ is the probability measure defined on $E_k$ such that, for any measurable function $f$ on $E_k$
$$
\eta_k(f):=\frac{\E\Big(f(X_k)\prod_{p=0}^{k-1}G_p(X_p)\Big)}{\E\Big(\prod_{p=0}^{k-1}G_p(X_p)\Big)}\ .
$$
%Then equation \ref{chg} is immediate by induction.
%
The measures $(\eta_k)_{0\leq k\leq n}$ defined above can be seen as the laws of $(X_k)_{0\leq k\leq n}$ under probability $\QQ$. Loosely speaking, the process $(X_k)_{0\leq k\leq n}$ with distribution $(\eta_k)_{0\leq k\leq n}$ is designed under the constrain $(\prod^k_{p=0}G_p)_{0\leq k\leq n}$. An intuitive interpretation comes by setting $G_k(x_k)=1_{A_k}(x_k)$ with $A_k\subset E_k$, then the process with distribution $\eta_k$ is just the ones surviving in the subsets $A_k$. It follows that the measures $\eta_k$ seem to be a relevant choice for the change of probability in our path dependent situation.\\
Furthermore, it is also important to observe that, for any measurable function $f$ on $E_k$
\begin{equation}
\label{eq:etaf}
\eta_k(f)=\frac{\eta_{k-1}(G_{k-1}M_k(f))}{\eta_{k-1}(G_{k-1})}\ .
\end{equation}
%or equivalently, for all $x_k\in E_k$, 
%$$
%\eta_k(dx_k)=\frac{\eta_{k-1}(G_{k-1}M_k(\cdot,dx_k))}{\eta_{k-1}(G_{k-1})}\ .
%$$
We denote the recursive relation between $\eta_k$ and $\eta_{k-1}$ by introducing the operators $\Phi_k$ such that, for all $1\leq  k\leq n$
\begin{equation}
\label{eq:Phi}
\eta_k=\Phi_k(\eta_{k-1})\ .
\end{equation}
Let us now introduce the integral operator $Q_{k}$ such that, for all $1\leq  k\leq n$
\begin{equation}
\label{eq:Q}
Q_{k}(f)(x_{k-1}):=\int G_{k-1}(x_{k-1})M_{k}(x_{k-1},dx_{k})f(x_{k})\ .
%\quad\textrm{or equivalently}\quad Q _{k}(x_{k-1},dx_{k})=G_{k-1}(x_{k-1})M_{k}(x_{k-1},dx_{k})\ .
\end{equation}
In further developments of this article, we suppose that $M_k(x_{k-1}, \cdot)$ are equivalent to some measures $\lambda_k$, for any $0\leq k\leq n$ and $x_{k-1}\in E_{k-1}$, i.e. there exists a collection of positive functions $H_k$ and measures $\lambda _k $ such that:
\begin{eqnarray}\label{hypo}
M_k(x_{k-1}, dx_k)=H_k(x_{k-1}, x_k)\lambda _k(dx_k)\ .
\end{eqnarray}
Now, we are in a position to state the following Lemma. 
\begin{lem}
\label{lem:RSE}
  For any measure $\eta$ on $E_k$, recursion~(\ref{RSE0}) defining $v_k$ can be rewritten:
%With this measure change, the recursion \ref{RSE0} turns into:
%Replacing the new expression \ref{expm} of $M_k(x_{k-1},dx_k)$ in the recursion \ref{RSE0}, we get:
\begin{equation*}%\label{RSE}
%v_k(x_k)=f_k(x_k)\vee\int\Phi_{k+1}(\eta)(dx_{k+1})\frac{G_k(x_k)H_{k+1}(x_k,x_{k+1})\eta(G_k)}{\eta(G_kH_{k+1}(\cdot,x_{k+1}))}v_{k+1}(x_{k+1})\ ,
v_k(x_k)=f_k(x_k)\vee Q_{k+1}(v_{k+1})(x_k)=f_k(x_k)\vee\Phi_{k+1}(\eta)\left (\frac{dQ _{k+1}(x_k,\cdot)}{d\Phi_{k+1}(\eta)}v_{k+1}\right )\ ,
\end{equation*}
for any $x_k\in E_k$,  where
$$
\frac{dQ _{k+1}(x_k,\cdot)}{d\Phi_{k+1}(\eta)}(x_{k+1})=\frac{G_k(x_k)H_{k+1}(x_k,x_{k+1})\eta(G_k)}{\eta(G_kH_{k+1}(\cdot,x_{k+1}))}\ ,
$$
for any $(x_k,x_{k+1})\in E_k\times E_{k+1}$. 
\end{lem}
\proof
Under Assumption~(\ref{hypo}), we have immediately the following formula
\begin{eqnarray}\label{expm}
M_{k+1}(x_{k},dx_{k+1})=H_{k+1}(x_{k},x_{k+1})\frac{\eta_{k}(G_{k})}{\eta_{k}(G_{k}H_{k+1}(\cdot,x_{k+1}))}\eta_{k+1}(dx_{k+1})\ .
\end{eqnarray}
Now, note that the above equation is still valid for any measure $\eta$, 
\begin{equation}\label{equ1}
M_{k+1}(x_{k},dx_{k+1})=H_{k+1}(x_{k},x_{k+1})\frac{\eta(G_{k})}{\eta(G_{k}H_{k+1}(\cdot,x_{k+1}))}\Phi_{k+1}(\eta)(dx_{k+1})\ .
\end{equation} 
Hence, the Radon Nikodym derivative of $M_{k+1}(x_{k},dx_{k+1})$ w.r.t. $\Phi_{k+1}(\eta)$ is such that
\begin{equation}
\frac{dM_{k+1}(x_{k},\cdot)}{d\Phi_{k+1}(\eta)}(x_{k+1})=H_{k+1}(x_{k},x_{k+1})\frac{\eta(G_{k})}{\eta(G_{k}H_{k+1}(\cdot,x_{k+1}))}\ .
\end{equation} 
We end the proof by applying above arguments to recursion~(\ref{RSE0}).

\section{A particle approximation scheme}\label{desmodel}
%\section{Description of the model}\label{desmodel}
%
From the above discussion,  we conclude that the distributions $(\eta_k)_{0\leq k\leq n}$ are a very good choice for the change of probability for the stochastic mesh model. In this section,  we first propose a particle model to sample the random variables according to these distributions, then we describe the resulting particle scheme proposed to approximate the Snell envelope $(v_k)_{0\leq k\leq n}$. \\
%In the precedent section, we have introduced the recursion $\eta_{k+1}=\Phi_{k+1}(\eta_k)$. 
By definition~(\ref{eq:Phi}) of $\Phi_{k+1}$, we have the following formula
\begin{eqnarray}
\Phi_k(\eta_{k-1})=\eta_{k-1}K_{k,\eta_{k-1}}=\eta_{k-1}S_{k-1,\eta_{k-1}}M_k=\Psi_{G_{k-1}}(\eta_{k-1})M_k\ .
\end{eqnarray} 
Where $K_{k,\eta_{k-1}}$, $S_{k-1,\eta_{k-1}}$ and $\Psi_{G_{k-1}}$ are defined as follows:
$$
\left \{
\begin{array}{lll}
%K_{k,\eta_{k-1}}(x,dz)=(S_{k-1,\eta_{k-1}}M_k)(x,dz)=\int S_{k-1,\eta_{k-1}}(x,dy)M_k(y,dz)\ .
K_{k,\eta_{k-1}}(x_{k-1},dx_k)&=&(S_{k-1,\eta_{k-1}}M_k)(x_{k-1},dx_k)\\\\
&=&\int S_{k-1,\eta_{k-1}}(x_{k-1},dx'_{k-1})M_k(x'_{k-1},dx_k)\ ,\\\\
S_{k-1,\eta_{k-1}}(x,dx')&=&\epsilon G_{k-1}(x)\delta_x(dx')+(1-\epsilon G_{k-1}(x))\Psi_{G_{k-1}}(\eta_{k-1})(dx')\\\\
\Psi_{G_{k-1}}(\eta_{k-1})(dx)&=&\frac{G_{k-1}(x)}{\eta_{k-1}(G_{k-1})}\eta_{k-1}(dx)\ ,
\end{array}
\right .
$$
where the real  $\epsilon$ is such that  $\epsilon G$ takes its values $[0,1]$. \\
%\begin{eqnarray*}
%&&S_{k-1,\eta_{k-1}}(x,dy)=\epsilon G_{k-1}(x)\delta_x(dy)+(1-\epsilon G_{k-1}(x))\underbrace{ \Psi_{G_{k-1}}(\eta_{k-1})(dy)}\\
%&&\hspace{101mm}||\\
%&&\hspace{90mm}=\frac{G_{k-1}(y)}{\eta_{k-1}(G_{k-1})}\eta_{k-1}(dy)\ .
%\end{eqnarray*}
More generally, the operations $\Psi$ and $S$ can be expressed as $\Psi_{G}(\eta)(f)=\frac{\eta(Gf)}{\eta(G)}=\eta S_{\eta}(f)$ with $S_{\eta}(f)=\epsilon Gf+(1-\epsilon G)\Psi_G(\eta)(f)$. \\
The particle approximation provided in the present paper is defined in terms of a Markov chain $\xi_k^{(N)}=(\xi^{(i,N)}_k)_{1\leq i\leq N}$ on the product state spaces $E^{N}_k$,
where the given integer $N$ is the number of particles sampled in every instant. 
The initial particle system, $\xi^{(N)}_0=\left(\xi^{(i,N)}_0\right)_{1\leq i\leq N}$, is a collection of  $N$ i.i.d. random copies of $X_0$. We let $\Fa^N_k$ be the sigma-field generated by the particle approximation model
from the origin, up to time $k$.
To simplify the presentation, when there is no confusion we suppress the population size parameter $N$, and we write $\xi_k$ and $\xi_k^{i}$
instead of $\xi_k^{(N)}$ and $\xi_k^{(i,N)}$. By construction, $\xi_k$ is a particle model with a selection transition
and a mutation type exploration i.e. the evolution from $\xi_k$ to $\xi_{k+1}$ is composed by two steps:
\begin{equation}\label{genetic transition}
\xi_{k}\in E_k^{N}~\overset{\rm Selection} {\underset{S_{k,\eta^N_{k}}} {-\!\!\!\!-\!\!\!\! 
-\!\!\!\!-\!\!\!\! 
-\!\!\!\!-\!\!\!\!-\!\!\!\!\longrightarrow}}~\widehat{\xi}_k:=\left(\widehat{\xi}^i_k\right)_{1\leq i\leq {N}}\in E_{k}^{N}  
~\overset{\rm Mutation}{\underset{M_{k+1}} {-\!\!\!\!-\!\!\!\!-\!\!\!\! 
-\!\!\!\!-\!\!\!\!-\!\!\!\!\longrightarrow}}~ 
\xi_{k+1}\in E_{k+1}^{N} \ .
\end{equation}
Then we define $\eta^N_k$ and $\widehat{\eta}^N_k$ as the occupation measures after the mutation and the selection steps. More precisely, 
$$
\eta^N_k:=\frac{1}{N}\sum_{1\leq i\leq N}\delta_{\xi^i_k}\quad\mbox{\rm and}\quad
\widehat{\eta}^N_k:=\frac{1}{N}\sum_{1\leq i\leq N}\delta_{\widehat{\xi}^i_k}\ .
$$
During the selection transition $S_{k,\eta^N_{k}}$, for $0\leq i\leq N$ with a probability $\epsilon G_k(\xi^i_k)$ we decide to skip the selection step i.e. we leave $\widehat{\xi}^i_k$ stay on particle $\xi^i_k$, and with probability $1-\epsilon G_k(\xi^i_k)$ we decide to do the following selection:  $\widehat{\xi}^i_k$ randomly takes the value in $\xi^j_k$ for $0\leq j\leq N$ with distribution  $\frac{G_k(\xi^j_k)}{\sum_{l=1}^NG_k(\xi^l_k)}$. Note that when $\epsilon G_k\equiv 1$,  the selection is skipped ( i.e. $\widehat{\xi}_k=\xi_k$) so that the model corresponds exactly to the Broadie-Glasserman type model analysed by P. Del Moral and P. Hu et al.~\cite{robust}. Hence, the factor $\epsilon$ can be interpreted as a level of selection against the rare events. \\
During the mutation transition $\widehat{\xi}_k\leadsto \xi_{k+1}$, every selected individual $\widehat{\xi}^i_k$ evolves randomly to a new individual $\xi_{k+1}^i=x$ randomly chosen with the distribution $M_{k+1}(\widehat{\xi}^i_k,dx)$, for $1\leq i\leq N$.\\
It is important to observe that by construction, $\eta^N_{k+1}$ is the empirical measure associated with $N$ conditionally independent and identically distributed random individual $\xi^i_{k+1}$ with common distribution $\Phi_{k+1}(\eta_k^N)$.

Now, we are in a position to describe precisely the new approximation scheme proposed to estimate the Snell envelope $(v_k)_{0\leq k\leq n}$. 
%If we denote by $Q _{k+1}(x_k,dx_{k+1})=G_k(x_k)M_{k+1}(x_k,dx_{k+1})$, and set $\eta=\eta_k^N$ in Proposition~\ref{RSE}, then 
The main idea consists in taking $\eta=\eta_k^N$, in Lemma~\ref{lem:RSE}, then observing that Snell envelope  $(v_k)_{0\leq k\leq n}$ is solution of the following recursion, for all $0\leq k<n$, 
%\begin{eqnarray}
%v_k(x_k)=f_k(x_k)\vee\int\Phi_{k+1}(\eta^N_k)(dx_{k+1}) \frac{dQ _{k+1}(x_k,\cdot)}{d\Phi_{k+1}(\eta^N_k)}(x_{k+1})v_{k+1}(x_{k+1})\ .
$$
v_k(x_k)=f_k(x_k)\vee\Phi_{k+1}(\eta^N_k)\left (\frac{dQ _{k+1}(x_k,\cdot)}{d\Phi_{k+1}(\eta^N_k)}v_{k+1}\right )\ .
$$
%\end{eqnarray} 
%
Now, if $\Phi_{k+1}(\eta^N_k)$ is well estimated by $\eta^N_{k+1}$, it is relevant to approximate $v_k$ by  $\widehat{v}_k$ defined by the following backward recursion
\begin{equation}
\label{eq:vEst}
\left \{
\begin{array}{lll}
\widehat{v}_n&=&f_n\\
\widehat{v}_k(x_k)&=&f_k(x_k)\vee  \eta^N_{k+1}\left ( {\displaystyle \frac{dQ_{k+1}(x_k,\cdot)}{d\Phi_{k+1}(\eta^N_k)} \widehat{v}_{k+1}} \right )\quad \textrm{for all} \ 0\leq k<n\ ,\\
%\widehat{v}_k(x_k)&=&f_k(x_k)\vee \int \eta^N_{k+1}(dx_{k+1}) \frac{dQ_{k+1}(x_k,\cdot)}{d\Phi_{k+1}(\eta^N_k)}(x_{k+1}) \widehat{v}_{k+1}(x_{k+1})\ \textrm{for all} \ 0\leq k<n\ ,\\
%&=&f_k(x_k)\vee\int\eta^N_{k+1}(dx_{k+1})\frac{G_k(x_k)H_{k+1}(x_k,x_{k+1})\eta^N_k(G_k)}{\eta^N_k(G_kH_{k+1}(\cdot,x_{k+1}))}\widehat{v}_{k+1}(x_{k+1})\ ,
\end{array}
\right .
\end{equation}
%for $0\leq k<n$, with terminal condition $\widehat{v}_n=f_n$. 
%where for all $(x_k,x_{k+1})\in E_k\times E_{k+1}$, 
%$$
%\frac{dQ_{k+1}(x_k,\cdot)}{d\Phi_{k+1}(\eta^N_k)}(x_{k+1}) =
%\frac{G_k(x_k)H_{k+1}(x_k,x_{k+1})\eta^N_k(G_k)}{\eta^N_k(G_kH_{k+1}(\cdot,x_{k+1}))}\ .
%$$
Note that in the above fomula~(\ref{eq:vEst}), the function  $v_k$ is defined not only on $E^N_k$ but on the whole state space $E_k$.\\
To simplify notations, we set 
$$
\widehat{Q}_{k+1}(x_k,dx_{k+1})=\eta^N_{k+1}(dx_{k+1})\frac{dQ_{k+1}(x_k,\cdot)}{d\Phi_{k+1}(\eta^N_k)}(x_{k+1})\ .
$$
Finally, with this notation, the real Snell envelope $(v_k)_{0\leq k\leq n}$ and the approximation $(\hat{v}_k)_{0\leq k\leq n}$ are such that, for all $0\leq k<n$, 
\begin{eqnarray*}
v_k&=&f_k\vee Q_{k+1}(v_{k+1})\\
\widehat{v}_k&=&f_k\vee \widehat{Q}_{k+1}(\widehat{v}_{k+1})\ .
\end{eqnarray*}

\section{Convergence and bias analysis}\label{convana}

By the previous construction, we can approximate $\Phi_{k+1}(\eta^N_k)$ by $\eta^N_{k+1}$. In this section, we will first analyze the  error associated with that approximation and then derive an error bound 
%of the error induced by
for the resulting Snell envelope approximation scheme. To simplify notations, in further development, we consider the random fields $V^N_k$ defined as 
$$
V^N_k:=\sqrt{N}~\left(\eta^N_k-\Phi_k(\eta^N_{k-1})\right)\ .
$$
The following lemma shows the conditional unbiasedness property and mean error estimates for the approximation $\eta^N_{k+1}$ of $\Phi_{k+1}(\eta^N_k)$. 
\begin{lem}~\label{lma}
For any integer $p\geq 1$, we denote by p' the smallest even integer greater than p. In this notation, for any $0\leq k\leq n$ and any integrable function $f$ on $E_{k+1}$, we have
\begin{eqnarray*}\label{mean-error}
\EE\left(\eta_{k+1}^N(f)|\Fa^N_k\right)=\Phi_{k+1}(\eta^N_k)(f)
\end{eqnarray*}
and
\begin{eqnarray*}
\EE\left(\left|V^N_k(f)\right|^p|\Fa^N_k\right)^{\frac{1}{p}}
\leq 2~a(p)~\left[\Phi_{k+1}(\eta^N_{k})(|f|^{p^{\prime}})\right]^{\frac{1}{p^{\prime}}}
\end{eqnarray*}
with the collection of  constants 
$$
a(2p)^{2p}=(2p)_p~2^{-p}\quad\mbox{and}\quad
a(2p+1)^{2p+1}=\frac{(2p+1)_{p+1}}{\sqrt{p+1/2}}~2^{-(p+1/2)}\ .
$$
\end{lem}
{\bf Proof :}
The conditional unbiasedness property is easily proved as follows 
\begin{eqnarray*}
\EE\left(\eta_{k+1}^N(f)|\eta^N_k\right)&=&\frac{1}{N}\sum^N_{i=1}\EE(f(\xi^i_{k+1})|\eta^N_k)\\
&=&\frac{1}{N}\sum^N_{i=1}K_{k+1,\eta^N_{k}}(f)(\xi^i_k)\\
&=&(\eta^N_kK_{k+1,\eta^N_k})(f)=\Phi_{k+1}(\eta^N_k)(f)\ .
\end{eqnarray*}
Then the above equality implies 
\begin{eqnarray*}
\EE\left(\left|\left[\eta_{k+1}^N-\Phi_{k+1}(\eta^N_k)\right](f)\right|^p|\Fa^N_k\right)^{\frac{1}{p}}\leq \EE\left(\left|\left[\eta_{k+1}^N-\mu_{k+1}^N\right](f)\right|^p|\Fa^N_k\right)^{\frac{1}{p}}\ ,
\end{eqnarray*}
where $\mu^{N}_{k+1}:=\frac{1}{N}\sum_{i=1}^{N}\delta_{Y^i_{k+1}}$ stands for an independent copy of $\eta_{k+1}^N$ given $\eta^N_{k}$.
Using Khintchine's type inequalities yields that
\begin{eqnarray*}
\sqrt{N}~\E\left( \left| [\eta_{k+1}^N-\mu^{N}_{k+1}](f) \right|^p
\left|\Fa^N_{k}\right.\right)^{\frac{1}{p}}
&\leq& 2~a(p)~\E\left(\left|f\left(\xi^1_{k+1}\right)\right|^{p^{\prime}}|~\Fa^N_{k}\right)^{\frac{1}{p^{\prime}}}\\&=&2~a(p)~\left[\Phi_{k+1}(\eta^N_{k})(|f|^{p^{\prime}})\right]^{\frac{1}{p^{\prime}}}\ .
\end{eqnarray*}
We end the proof by combining the above two inequalities.

\cqfd
A consequence of the unbiasedness property proved in Lemma~\ref{lma} is that
\begin{eqnarray*}
\E(\widehat{Q}_{k+1}(f)(x_k)|\eta^N_k)=Q_{k+1}(f)(x_k)\ .
\end{eqnarray*}
To estimate the error between $v_k$ and the approximation $\hat{v}_k$, it is usefull to introduce the following random integral operator $R^N_k$ such that for any measurable function on $E_{k+1}$, 
$$
R^{N}_{k+1}(f)(x_k)=\sqrt{N}\left(\widehat{Q}_{k+1}(f)(x_k)-Q_{k+1}(f)(x_k)\right)\ .
$$
Note that
$$
R^{N}_{k+1}(f)(x_k):=\int V^N_{k+1}(dx_{k+1})~\frac{dQ_{k+1}(x_k,\point)}{d\Phi_{k+1}(\eta^N_k)}(x_{k+1})~f(x_{k+1})\ ,
$$
then, applying again Lemma~\ref{lma} implies the following  Khintchine's type inequality 
\begin{eqnarray*}
&&\EE(\left|R^{N}_{k+1}(v_{k+1})(x_k)\right|^p|\eta^N_k)^{\frac{1}{p}}\\
&&\hspace{3mm}\leq 2~a(p)\left[\int_{E_{k+1}}\Phi_{k+1}(\eta^N_k)(dx_{k+1})\left(\frac{dQ_{k+1}(x_k,\cdot)} {d\Phi_{k+1}(\eta^N_k)}(x_{k+1}) v_{k+1}(x_{k+1})\right)^{p'}\right]^{\frac{1}{p'}}
\end{eqnarray*}
Let $\widehat{Q}_{k,l}=\widehat{Q}_{k+1}\widehat{Q}_{k+2}\ldots \widehat{Q}_{l}$ for any $0\leq k<l\leq n$,  then it follows easily, by recursion, that  
$$
\E(\widehat{Q}_{k,l}(f)(x_k)|\eta^N_k)=Q_{k,l}(f)(x_k)\ .
$$
Now, by Lemma~\ref{lem:lemme1}, we conclude 
\begin{equation}\label{decom}
\sqrt{N}\left|(v_k-\widehat{v}_k)\right|\leq \sum_{k<l<n}\widehat{Q}_{k,l}|(R^{N}_{l+1})(v_{l+1})|\ .
\end{equation}
We are now in position to state the main result of this paper. 
\begin{thm}\label{thethm} 
For any $0\leq k\leq n$ and any integer $p\geq 1$, we have
\begin{eqnarray*}
\sup_{x\in E_k}\left\Vert (\widehat{v}_k-v_k)(x)\right\Vert _{L_p}\leq\sum_{k<l<n}\frac{2~a(p)} {\sqrt{N}}q_{k,l}  ~\left[Q_{k,l+1}(h_{l+1}^{p'-1}v_{l+1}^{p'})(x)\right]^{\frac{1}{p'}}\ ,
\end{eqnarray*}
with a collection of constants $q_{k,l}$ and functions $h_k$ defined as 
\begin{equation}
\label{eq:qh}
q_{k,l}:=\left[\Vert G_{l}\Vert ~\Vert h_{k+1}\Vert \prod^{l-1}_{m=k}\Vert G_{m}\Vert \right]^{\frac{p'-1}{p'}}\quad\textrm{and}\quad h_k(x_k):=\sup_{x,y\in E_{k-1}}\frac{H_{k}(x,x_{k})}{H_{k}(y,x_{k})}\ .
\end{equation}
\end{thm}
{\bf Proof :}
First, decomposition~(\ref{decom}) yields 
$$
\sqrt{N}\left\Vert (\widehat{v}_k-v_k)(x)\right\Vert_{L_p}\leq\sum_{k<l<n}\left \Vert\widehat{Q}_{k,l}\vert (R^{N}_{l+1})(v_{l+1})\vert (x)\right\Vert_{L_p}\ ,\quad \textrm{for all}\ x\in E_k\ .
$$
Note that 
$$
\Vert \widehat{Q}_{k,l}(1)\Vert \leq b_{k,l}\ ,\quad\textrm{where}\quad b_{k,l}:=\Vert h_{k+1}\Vert \prod^{l-1}_{m=k}\Vert G_{m}\Vert \ .
$$
Then it follows easily that for any integrable function $f$ on $E_l$
$$
(\widehat{Q}_{k,l}(f))^p\leq(b_{k,l})^{p-1}\widehat{Q}_{k,l}(f^p)\ .
$$
This yields that
$$
\left \Vert \widehat{Q}_{k,l}\left|(R^{N}_{l+1}))(v_{l+1})\right|(x)\right\Vert_{L_p}\leq(b_{k,l})^{\frac{p-1}{p}}\EE\left(\widehat{Q}_{k,l}\left(\left|(R^{N}_{l+1}))(v_{l+1})\right|\right)^p(x)\right)^{\frac{1}{p}}\ .
$$
Applying Lemma~\ref{lma} to the right-hand side of the above inequality, we obtain for any $x_l\in E_l$
\begin{eqnarray*}
&&\EE\left(\left|(R^{N}_{l+1}))(v_{l+1})(x_l)\right|^p|\eta^N_l\right)^{\frac{1}{p}}\\
&&\hspace{3mm}\leq 2~a(p)\left[\int_{E_{l+1}}\Phi_{l+1}(\eta^N_l)(dx_{l+1})\left(\frac{dQ_{l+1}(x_l,\cdot)} {d\Phi_{l+1}(\eta^N_l)}(x_{l+1}) v_{l+1}(x_{l+1})\right)^{p'}\right]^{\frac{1}{p'}}
\end{eqnarray*}
from which we find that 
\begin{eqnarray*}
&&\EE\left(\left|(R^{N}_{l+1}))(v_{l+1})(x_l)\right|^p|\eta^N_l\right)^{\frac{1}{p}}\\
&&\hspace{3mm}\leq 2~a(p)\left[\int_{E_{l+1}}Q_{l+1}(x_l, dx_{l+1})\left(\frac{dQ_{l+1}(x_l,\cdot)} {d\Phi_{l+1}(\eta^N_l)}(x_{l+1}) \right)^{p'-1}v_{l+1}(x_{l+1})^{p'}\right]^{\frac{1}{p'}}
\end{eqnarray*}
By definition~(\ref{eq:qh}) of functions $h_{l+1}$ and in developing the Radon Nikodym derivative, we obtain 
$$
\frac{dQ_{l+1}(x_l,\cdot)} {d\Phi_{l+1}(\eta^N_l)}(x_{l+1})=\frac{\eta^N_l(G_l)G_l(x_l)H_{l+1}(x_l,x_{l+1})}{\eta^N_l(G_lH_{l+1})(\cdot,x_{l+1})}\leq \Vert G_l\Vert h_{l+1}(x_{l+1})\ ,
$$
which implies 
\begin{eqnarray*}
&&\EE\left(\left|(R^{N}_{l+1}))(v_{l+1})(x_l)\right|^p|\eta^N_l\right)^{\frac{1}{p}}\\
&&\hspace{3mm}\leq 2~a(p)\Vert G_l\Vert ^{\frac{p'-1}{p'}}\left[\int_{E_{l+1}}Q_{l+1}(x_l, dx_{l+1})\left(h_{l+1}(x_{l+1})\right)^{p'-1}v_{l+1}(x_{l+1})^{p'}\right]^{\frac{1}{p'}}
\end{eqnarray*}
Gathering the above arguments, we conclude that
\begin{eqnarray*}
\left \Vert \left(\widehat{v}_k-v_k\right)(x)\right\Vert_{L_p}
\leq\sum_{k<l<n}\frac{2~a(p)} {\sqrt{N}}q_{k,l}  ~\left(Q_{k,l+1}(h_{l+1}^{p'-1}v_{l+1}^{p'})(x)\right)^{\frac{1}{p'}}\ .
\end{eqnarray*}

\cqfd
{\bf Remarks :}
The constants $q_{k,l}$ could be largely reduced. In fact, $q_{k,l}$ comes from bounding $\Vert \prod_m\eta^N_m(G_m)\Vert _{L_p}$. In~\cite{l2}, the authors proved $\Vert \prod_mG_m\Vert _{L_2}+\frac{constant}{N}$ as a non asymptotic boundary for $\Vert \prod_m\eta^N_m(G_m)\Vert _{L_2}$. In most cases, the functions $G$ take their values in $[0,1]$, then the majoration $\Vert \prod_mG_m\Vert \leq 1$ holds, but $\Vert \prod_mG_m\Vert _{L_2}$ is very small.

When the function $G$ vanishes in some regions of the state space, we also mention that the particle model is only defined up to the first time $\tau^N = k$ such that $\eta_k^N(G_k) = 0$. We can prove that the event $\{\tau ^N \leq n\}$ has an exponentially small probability to occur, with the number of particles $N$. In fact, the estimates presented in the above theorems can be extended to this singular situation by replacing $\widehat{v}_k$ by the particle estimates $\widehat{v}_k1_{\tau^N \geq n}$. The stochastic analysis of these singular models are quite technical, for further details we refer the reader to section 7.2.2 and section 7.4 in the book~\cite{fk}.

To understand better the $\LL_p$-mean error bounds in the above theorem, we deduce the following exponential concentration inequality

\begin{prop}
For any $0\leq k\leq n$ any and any $\epsilon>0$, we have
\begin{eqnarray}\label{expon}
\sup_{x\in E_k}{\PP\left(|v_k(x)-\widehat{v}_k(x)|>\frac{c}{\sqrt{N}}+\epsilon\right)}\leq exp\left({-N\epsilon^2}/{c^2)} \right)\ ,
\end{eqnarray}
with constant $c=\sum_{k<l<n}2~ q_{k,l} ~\left(Q_{k,l+1}(h_{l+1}^{p'-1}v_{l+1}^{p'})(x)\right)^{\frac{1}{p'}}$. 
\end{prop}
{\bf Proof :}
This result is a direct consequence from the fact that for any non negative random variable $U$ such that
$$
\exists b<\infty~\mbox{\rm s.t.}~\forall r\geq 1\qquad \EE\left(U^r\right)^{\frac{1}{r}}\leq  a(r)~b~\Rightarrow
\PP\left(U\geq b+\epsilon\right)\leq \exp{\left(-{\epsilon^2}/{(2b^2)}
                                                                           \right)}\ . 
$$
To check this claim, we develop the exponential and verify that
\begin{eqnarray*}
\forall t\geq 0~~ \EE\left(e^{tU}\right)
\leq \exp{\left(\frac{(bt)^2}{2}+bt\right)} \Rightarrow\PP(U\geq b+\epsilon)\leq \exp{\left(-\sup_{t\geq 0}{(\epsilon t-\frac{(bt)^2}{2})}\right)}
\end{eqnarray*}
\cqfd

Simarly to the orginal Broadie-Glasserman model, the following proposition shows that in this model we also over-estimate the Snell envelope.
\begin{prop}
For any $0\leq k\leq n$ and any $x_k\in E_k$
\begin{equation}
\EE\left(\widehat{v}_k(x_k)\right)\geq v_k(x_k)\ .
\end{equation}
\end{prop} 
\proof
We can easily prove this inequality with a simple backward induction. The terminal condition $\widehat{v}_n=v_n$ implies directly the inequality at instant $n$. Assuming the inequality at time $k+1$, then the Jensen's inequality implies
\begin{eqnarray*}
\EE\left(\widehat{v}_k(x_k)\right)&\geq& f_k(x_k)\vee \EE\left(\widehat{Q}_{k+1}\widehat{v}_{k+1}(x_k)\right) \\
&=& f_k(x_k)\vee \EE\left(\int_{E^N_{k+1}}\widehat{Q}_{k+1}(x_k, dx_{k+1})\EE\left(\widehat{v}_{k+1}(x_{k+1})|\Fa^N_{k+1}\right)\right)\ .
\end{eqnarray*}
By the induction assumption at time $k+1$, we have
\begin{eqnarray*}
\EE\left(\int_{E^N_{k+1}}\widehat{Q}_{k+1}(x_k, dx_{k+1})\EE\left(\widehat{v}_{k+1}(x_{k+1})|\Fa^N_{k+1}\right)\right)&\geq&\EE\left(\widehat{Q}_{k+1}v_{k+1}(x_k)\right)\\
&=& Q_{k+1}v_{k+1}(x_k)\ .
\end{eqnarray*}
Then the inequatily still holds at time $k$, which completes the proof.

\cqfd

\end{document}